\documentclass[a4paper]{article}
\usepackage{amsfonts}
\usepackage{amssymb}
\usepackage{latexsym}
\newtheorem{theorem}{Theorem}[section]

\newtheorem{lemma}[theorem]{Lemma}
\newtheorem{definition}[theorem]{Definition}
\newtheorem{remark}[theorem]{Remark}

\newcommand{\proof}{\medskip \noindent {\bf Proof. \ \ }}
\newcommand{\qed}{\null\hfill $\Box\;\;$ \medskip}

\begin{document}

\parbox{1mm}

\begin{center}
{\bf {\sc \Large The Egoroff Theorem for Operator-Valued Measures
in Locally Convex Spaces}}
\end{center}

\vskip 12pt

\begin{center}
{\bf J\'an Halu\v{s}ka and Ondrej Hutn\'ik}\footnote{{\it
Mathematics Subject Classification (2000):} Primary 46G10;
Secondary 06F20
\newline {\it Key words and phrases:} Operator valued measure,
locally convex topological vector spaces, Egoroff theorem,
convergence in measure, net convergence of functions.
\newline {\it Acknowledgement.} This paper was partially supported by Grants
VEGA 2/0097/08 and VVGS 45/10-11.}
\end{center}

\vskip 24pt

\hspace{5mm}\parbox[t]{11cm}{\fontsize{9pt}{0.1in}\selectfont\noindent{\bf
Abstract.} The Egoroff theorem for measurable $\bold X$-valued
functions and operator-valued measures $\bold m: \Sigma \to
L(\bold X, \bold Y)$, where $\Sigma$ is a $\sigma$-algebra of
subsets of $T \neq \emptyset$ and $\bold X$, $\bold Y$ are both
locally convex spaces, is proved. The measure is supposed to be
atomic and the convergence of functions is net.} \vskip 24pt

\section{Introduction}

The classical Egoroff theorem states that almost everywhere
convergent sequences of measurable functions on a finite measure
space converge almost uniformly, that is, for every $\varepsilon >
0$ the convergence is uniform on a set whose complement has
measure less than $\varepsilon$, cf.~\cite{Halmos}. A necessary
and sufficient condition for a sequence of measurable real
functions to be almost uniformly convergent is given
in~\cite{Bartle2}. In generalizing to functions taking values in
more general spaces, some problems appear arising from the fact
that the classical relationship between the pointwise convergence
and the convergence in measure is not saved. Namely, any
convergence almost everywhere does not imply convergence in
measure in general. For a complete measure space $(T, \Sigma,
\mu)$, and a locally convex space $\bold X$ the following
condition may be considered, cf.~\cite{Goguadze}: 

{\sl Let $\mathcal{M}$ be a family of $\bold X$-valued functions
defined on $T$. The locally convex space $\bold X$ is said to
satisfy the finite Egoroff condition with respect to $\mathcal{M}$
if and only if every sequence in $\mathcal{M}$ which converges
almost everywhere to a function, is almost uniformly convergent on
every measurable set $A \in \Sigma$ of a finite measure.}

However, it is not so easy to find such a family of functions
$\mathcal{M}$. Thus, E.~Wagner and W.~Wilczy\'nski proved the
following theorem in~\cite{WilWag}.

\begin{theorem}
If a measurable space $(T, \Sigma, \mu)$ fulfils the countable
chain condition, then the convergence $I$-a.e. is equivalent to
the convergence with respect to the $\sigma$-ideal $I$ if and only
if $\Sigma/I$ is atomic.
\end{theorem}

In the operator-valued measure theory in Banach spaces the
pointwise convergence (of sequences) of measurable functions on a
set of finite semivariation implies the convergence in
(continuous) semi\-va\-ria\-tion of the measure $\bold{m}:\Sigma
\to L(\bold X, \bold Y)$, where $\Sigma$ is a $\sigma$-algebra of
subsets of a set $T \neq \emptyset$, and $\bold X, \bold Y$  are
Banach spaces, cf.~\cite{Bartle}, and~\cite{Dobrakov}. If $\bold
X$ fails to be metrizable, the relationship between these two
convergences is quite unlike the classical situation, cf.~Example
after Definition~1.11 in~\cite{Smith-Tucker}.

It is also well-known that the Egoroff theorem cannot hold for
arbitrary nets of measurable functions without some restrictions
on measure, net convergence of functions, or class of measurable
functions. For instance, the net of characteristic functions of
finite subsets of $[0, 1]$ shows that one needs very special
measures (supported by a countable set) to have almost uniform
convergence. In~\cite{Hal91}, Definition~1.2, the first author
introduced the so called Condition~(GB) under which everywhere
convergence of net of measurable functions implies convergence of
these functions in semi\-va\-ria\-tion on a set of finite
va\-ria\-tion of measure in locally convex setting,
cf.~\cite{Hal94}, Theorem~3.3 (without the assumption of countable
chain condition). This condition concerns families of submeasures
and enables to work with nets of measurable functions instead of
sequences. The Condition~(GB) is a generalization of Condition~(B)
which is important for sequences in the classical measure and
integration theory, cf.~\cite{Luzin}. The analogous condition for
nets in the classical setting was introduced and investigated by
B.~F. Goguadze, cf.~\cite{Goguadze}.

Recall that Condition~(GB) is fulfilled in the case of atomic
operator-valued measures, cf.~\cite{Hal94}. Atomic measures are
not of a great interest in the classical theory of measure and
integral, because they lead only to considerations of absolutely
convergent series. But when we consider measures with very general
range space, e.g. a locally convex space, the situation changes.

In this paper we prove the Egoroff theorem for atomic
operator-valued measures in locally convex topological vector
spaces.

\section{Preliminaries}

The description of the theory of locally convex topological vector
spaces may be found in~\cite{Jar81}. By a {\it net} (with values
in a set $S$) we mean a function from $I$ to $S$, where $I$ is a
directed partially ordered set. Throughout this paper $I$ is a
directed index set representing direction of a net. For
terminology concerning nets, see~\cite{Keley}.

Let $T$ be a non-void set and let $\Sigma$ be a $\sigma$-\-algebra
of subsets of $T$. By $2^T$ we denote the potential set of $T$.
Let $\bold X$ and $\bold Y$ be two Hausdorff locally convex
topological vector spaces over the field $\mathbb{K}$ of all real
$\mathbb{R}$ or complex numbers $\mathbb{C}$, with two families of
seminorms $P$ and $Q$ defining the topologies on $\bold X$ and
$\bold Y$, respectively. Let $L(\bold X, \bold Y)$ denote the
space of all continuous linear operators $L: \bold X\to \bold Y$,
and let $\mathbb{N}$ be the set of all natural numbers.

In what follows $\bold m: \Sigma \to L(\bold X, \bold Y)$ is an
ope\-ra\-tor-valued measure $\sigma$-ad\-di\-ti\-ve in the strong
operator topology of the space $L(\bold X, \bold Y)$, i.e. $\bold
m(\cdot) \bold x: \Sigma \to \bold Y$ is a $\bold Y$-valued vector
measure for every $\bold x \in \bold X$.

\begin{definition}\rm
Let $p \in P$, $q \in Q$, $E \in \Sigma$ and $\bold m: \Sigma \to
L(\bold X, \bold Y)$.
\begin{itemize}
\item[(a)] The $(p, q)$-{\it semi\-va\-ria\-tion of the measure}
$\bold m$ is the set function $\hat{\bold m}_{p, q}:\Sigma \to [0,
\infty ]$, defined as follows:
$$\hat{\bold m}_{p, q} (E) =  \sup q \left( \sum_{n = 1}^{N} \bold m(E_n \cap E) \bold x_n \right),$$
where the supremum is taken over all finite disjoint partitions
$\{E_n  \in \Sigma; \,E \subset \bigcup_{n=1}^{N} E_n, \,E_n  \cap
E_m = \emptyset, \,n \neq m,\, m, n = 1, 2, \dots , N \}$ of $E$,
and all finite sets $\{ \bold x_n \in \bold X; \,p(\bold x_n) \leq
1, \, n = 1, 2, \dots , N \}$, $N \in \mathbb{N}$. \item[(b)] The
$(p, q)$-{\it variation of the measure} $\bold m$ is the set
function $\bold{var}_{p, q} (\bold m, \cdot): \Sigma \to [0,
\infty ]$, defined by the equality
$$\bold{var}_{p, q}(\bold m, E) = \sup \sum_{n = 1}^{N} q_p \left(\bold m (E_n
\cap E)\right),$$ where the supremum is taken over all finite
disjoint partitions $\{E_n  \in \Sigma; \, E = \bigcup_{n=1}^{N}
E_n, \, E_n  \cap E_m  = \emptyset, \,n \neq m, \,n, m = 1, 2,
\dots , N, \, N \in \mathbb{N} \}$ of $E$ and $$q_p(\bold m (F)) =
\sup_{p(\bold x) \leq 1} q(\bold m(F)\bold x), \quad F \in
\Sigma.$$ \item[(c)] The {\it inner $(p,q)$-semivariation of the
measure} $\bold m$ is the set function $\hat{\bold m}^*_{p,q}: 2^T
\to [0, \infty ]$, defined as follows:
$$\hat{\bold m}^*_{p,q}(E) = \sup_{F \subset E, F \in \Sigma}
\hat{\bold m}_{p,q}(F), \quad E \in 2^T.$$
\end{itemize}
\end{definition}
Note that $\hat{\bold m}_{p, q}(E) \leq \bold{var}_{p, q}(\bold m,
E)$ for every $q \in Q$, $p \in P$, and $E \in \Sigma$. The
following lemma is obvious.

\begin{lemma}\label{lemma2.2}
For every $p \in P$ and $q \in Q$, the $(p, q)$-(semi)variation of
the measure $\bold m$ is a monotone and $\sigma$-additive
($\sigma$-subadditive) set function, and $\bold{var}_{p, q}(\bold
m, \emptyset)=0$, ($\hat{\bold m}_{p, q}(\emptyset)= 0)$.
\end{lemma}

\begin{definition}\rm Let $\bold m: \Sigma \to L(\bold X, \bold
Y)$.
\begin{itemize} \item[(a)] The set $E \in \Sigma$ is said
{\it to be of positive variation of the measure} $\bold m$ if
there exist $q \in Q$, $p \in P$, such that $\bold{var}_{p,
q}(\bold m, E) > 0$. 
\item[(b)] We say that the set $E \in \Sigma$ {\it is of finite
variation of the measure} $\bold m$ if to every $q \in Q$ there
exists $p \in P$, such that $\bold{var}_{p, q} (\bold m, E)<
+\infty$. We will denote this relation shortly $Q \to_E P,$ or,
$q\mapsto_E p$, for $q \in Q$, $p \in P$.
\end{itemize}
\end{definition}

\begin{remark}\rm
The relation $Q \to_E P$ may be different for different sets $E
\in \Sigma$ of finite variation of the measure $\bold m$.
\end{remark}

\begin{definition}\label{Definition 1.1} \rm
A measure $\bold m: \Sigma \to L(\bold X, \bold Y)$ is said to
{\it satisfy Condition~(GB)} if for every $E \in \Sigma$ of finite
and positive variation and every net $(E_i)_{i\in I}$, $E_i
\subset E$, of sets from $\Sigma$ there holds
$$\limsup_{i \in I} E_i \neq \emptyset,$$
whenever there exist real numbers $\delta (q, p, E) > 0$, $p \in
P$, $q \in Q$, such that $\hat {\bold m}_{p,q}(E_i) \geq \delta
(q, p, E)$ for every $i \in I$ and every couple $(p, q) \in P
\times Q$ with $q \mapsto_E p$.
\end{definition}

Note that the everywhere convergence of a net of measurable
functions implies the convergence of these functions in
semivariation on a set of finite variation if and only if
Condition~(GB) is fulfilled. For instance, the Lebesque measure
does not satisfy Condition~(GB) on the real line (for arbitrary
nets of sets).

\begin{definition}\rm
We say that a set $E \in \Sigma$ of positive semi\-va\-ria\-tion
of the measure $\bold m$ is an $\hat{\bold m}$-{\it atom} if every
proper subset $A$ of $E$ is either $\emptyset$ or $A \notin
\Sigma$.

We say that the measure $\bold m$ is {\it atomic} if each $E \in
\Sigma$ can be expressed in the form $E=\bigcup_{k = 1}^{\infty}
A_k$, where $A_k, k \in \mathbb{N}$, are $\hat{\bold m}$-atoms.
\end{definition}

In~\cite{Hal94} it is shown that a class of measures satisfying
Condition~(GB) is non-empty and the following result is proved.

\begin{theorem}\label{thmatomic}
If $\bold m$ is a (countable) purely atomic operator-valued
measure, then $\bold m$ satisfies Condition~(GB).
\end{theorem}

\begin{definition}\rm
We say that a function $\bold f: T \to \bold X$ is {\it
measurable} if $$\{t \in T; \,p(\bold f(t)) \geq  \eta\} \in
\Sigma$$ for every $\eta > 0$ and $p \in P$.
\end{definition}



\begin{definition}\label{Definition 1.2}\rm
A net $(\bold f_i)_{i \in I}$ of measurable functions is said to
be {\it $\bold m$-almost uniformly convergent} to a measurable
function $\bold f$ on $E \in \Sigma$ if for every $\varepsilon >
0$ and every $(p,q) \in P\times Q$ with $q\mapsto_{E} p$ there
exist measurable sets $F=E(\varepsilon, p, q)$, such that
$$\lim_{i \in I} \|\bold f_i - \bold f\|_{E \setminus
F,p} = 0\quad \textrm{and} \quad \hat{\bold m}_{p, q}(F) <
\varepsilon,$$ where $\|\bold g\|_{G,p}=\sup\limits_{t\in
G}p(\bold g(t))$.
\end{definition}

In~\cite{Hal97} the concept of generalized strong continuity of
the semi\-va\-ria\-tion of the measure is introduced. Note that
this notion enables development of the concept of an integral with
respect to the $L(\bold X, \bold Y)$-valued measure based on the
net convergence of simple functions. For this purpose the notion
of the inner semivariation is used for this generalization. This
way we restrict the set of $L(\bold X, \bold Y)$-valued measures
which can be taken for such type of integration. For instance,
every atomic measure is generalized strongly continuous. So, the
class of measures with the generalized strongly continuous
semivariation is nonempty.

\begin{definition}\rm\label{defGS}
We say that the semi\-va\-ria\-tion of the measure $\bold m:
\Sigma \to L(\bold X, \bold Y)$ is {\it generalized strongly
contin\-uous ({\rm GS}-continuous, for short)} if for every set of
finite variation $E \in \Sigma$ and every monotone net of sets
$(E_i)_{i\in I} \subset T$, $E_i \subset E$, $i \in I$, the
following equality  $$\lim_{i \in I} \hat{\bold m}_{p,q}^*(E_i) =
\hat{\bold m}_{p,q}^*\left(\lim_{i \in I} E_i\right)$$ holds for
every couple $(p, q) \in P \times Q $, such that $q \mapsto_E p$.
\end{definition}

The following result connecting the notion of GS-continuity and
Condition~(GB) was proved in~\cite{Hal97}.

\begin{theorem}\label{thmGB}
If the semivariation of a measure $\bold m: \Sigma \to L(\bold X,
\bold Y)$ is {\rm GS}-contin\-uous, then the measure $\bold m$
satisfies Condition~(GB).
\end{theorem}

\section{Egoroff theorem for atomic operator-valued measures in locally convex spaces}

Before proving the main result of this paper we state the
following useful lemma.

\begin{lemma}\label{lemma}
If $\bold m: \Sigma \to L(\bold X, \bold Y)$ is a (countable)
purely atomic measure, then its semivariation is {\rm
GS}-continuous.
\end{lemma}

\proof Let $E \in \Sigma$ be a set of finite and positive
variation of the measure $\bold m$. Let $(E_i)_{i \in I}$ be an
arbitrary decreasing net of sets from $\Sigma$. Recall that $E_i
\searrow G (\in 2^G)$ if and only if
\begin{itemize}
\item[(1)] $i \leq j \Rightarrow E_i \supset E_j$, and \item[(2)]
$\bigcap_{i \in I} E_i = G$.
\end{itemize}It is clear that it is enough to consider the case $G =
\emptyset$, because $$E_i \searrow G \Leftrightarrow G \subset
E_i, \quad E_i \setminus G \searrow \emptyset.$$ First, in the
case $E_i \in \Sigma$, $i \in I,$ we have
\begin{equation}\label{4}
\lim_{i \in I} \hat{\bold m}^*_{p, q}(E_i) = \lim_{i \in I}
\hat{\bold m}_{p, q}(E_i),
\end{equation}
and since the family of atoms is at most a countable set, there is
$$\lim_{i \in I} E_i = \bigcap_{i \in I} E_i \in \Sigma$$ and,
therefore,
\begin{equation}\label{5}
\hat{\bold m}^*_{p,q}\left(\lim_{i \in I} E_i\right) = \hat{\bold
m}_{p, q}\left(\lim_{i \in I} E_i\right)
\end{equation} for every $p \in P$, $q \in Q$ such that $q \to_{E} p$.

Now, take an arbitrary set $E \in \Sigma$ of positive and finite
variation. Denote by $\mathcal{A}$ the set of all $\hat{\bold
m}$-atoms, and put $l(i, E) = (\mathcal{A} \cap E) \setminus E_i$,
$i \in I$. Clearly $$i \leq j, \,\,i, j \in I \Rightarrow l(i, E)
\subset l(j, E)$$ and there exist atoms  $A_n \in \mathcal{A}$,
$n\in \mathbb{N}$, such that $l(i, E) = \{ A_1, A_2, \dots , A_n,
\dots\}$.

By Lemma~\ref{lemma2.2} we have
$$\bold{var}_{p,q}(\bold{m},E) = \bold{var}_{p,q}(\bold{m},E_{i}) +
\sum_{A_{n}\in l(i,E)} \bold{var}_{p,q}(\bold{m},A_{n}),$$ for
$i\in I$, and $p\in P$, $q\in Q$. Since
$$\hat{\bold{m}}_{p,q}(E_{i}) \leq \bold{\bold{var}}_{p,q}(\bold{m},E_{i}), \quad i\in I, p\in P, q\in
Q,$$ there is
$$\hat{\bold m}_{p,q}(E_i) \leq \bold{var}_{p,q}(\bold{m},E) -
\sum_{A_n \in l(i, E)} \bold{var}_{p, q}(\bold{m},A_n),$$ where
$i\in I$, $p \in P$, $q \in Q$. Since
$\bold{var}_{p,q}(\bold{m},E\cap\cdot): \Sigma \to [0,\infty)$ is
a finite real measure for every $p\in P$, $q\in Q$ with
$q\mapsto_{E}p$, then for every $\varepsilon >0$, $p\in P$, $q\in
Q$ such that $q\mapsto_{E}p$, there exists an index
$i_{0}=i_{0}(\varepsilon, p,q,E) \in I$, such that
\begin{equation}\label{6}
\hat{\bold{m}}_{p,q}(E_{i}) < \varepsilon
\end{equation} holds for every $i \geq i_{0}$, $i\in I$. Combining~(\ref{4}), (\ref{5}), (\ref{6})
and Definition~\ref{defGS} we see that the assertion is proved for
the case when $(E_i)_{i\in I}$ is a decreasing net of sets from
$\Sigma$. The other cases of monotone nets of sets may be proved
analogously.

Let now $G \subset T$ be an arbitrary set. Then there is exactly
one (countable) set $F^* = \mathcal{A} \cap G$ with the property:
$$\hat{\bold m}^*_{p, q}(G) =
\sup_{F \subset G, F \in \Sigma} \hat{\bold m}_{p, q}(F) =
\hat{\bold m}_{p, q}(F^*), \quad p \in P, q \in Q.$$

The proof for the inner measure and the arbitrary net of subsets
$(E_i)_{i \in I}$ goes by the same procedure as in the previous
part of proof concerning the set system $\Sigma$. \qed

\begin{remark}\rm
Observe that Lemma~\ref{lemma} and Theorem~\ref{thmGB} imply the
statement of Theorem~\ref{thmatomic}.
\end{remark}

\begin{theorem}[Egoroff]\label{Theorem 2.1}
Let $\bold m: \Sigma \to L(\bold X, \bold Y)$ be a purely atomic
measure, and let $E \in \Sigma$ be of finite variation of the
measure $\bold m$. Let $\bold f: T \to \bold X$ be a measurable
function, and $(\bold f_i: T \to \bold X)_{i \in I}$ be a net of
measurable functions, such that
\begin{equation}\label{net1}
\lim_{i \in I} p(\bold f_i(t) - \bold f(t)) = 0 \quad \textrm{for
every\,\,\,} t \in E\,\, \textrm{and\,\,\,} p \in P.
\end{equation}
Then a net $(\bold f_i)_{i \in I}$ of functions $\bold m$-almost
uniformly converges to $\bold f$ on $E \in \Sigma$.
\end{theorem}

\proof We have to prove that for a given $\varepsilon > 0$ and
every $q \in Q$, $p \in P$, such that $q \to_E p$ there exist
measurable sets $F = E(\varepsilon, p, q) \in \Sigma$, such that
$\lim_{i \in I} \| \bold f_i -   \bold f \|_{E \setminus F,p}=0,$
and $\hat{\bold m}_{p, q}(F) < \varepsilon$.

Suppose that~(\ref{net1}) holds. For every $m \in \mathbb{N}$, $p
\in P$, and $j \in I$, put {\setlength\arraycolsep{2pt}
\begin{eqnarray*}
B^p_{m, j} & = & E \cap \left\{ t \in T; \,p(\bold f_i(t)
- \bold f(t)) <  \frac{1}{m}, \quad i \geq j \right\} \\
& = & E \cap \bigcap_{i \geq j} \left\{t \in T; \,p(\bold f_i(t) -
\bold f(t)) < \frac{1}{m}, \quad i \in I \right\}.
\end{eqnarray*}}Since there are countable many of atoms, $B^p_{m, j} \subset
\Sigma$ and $\# B^p_{m, j} =  \aleph_0$. Clearly, if $i, j \in I$
such that $i \leq j$, then $ B^p_{m, i} \subset B^p_{m, j}$ for
every $m \in \mathbb{N}$ and $p \in P$. Put $$E^p_m = \bigcup_{j
\in I} B^p_{m, j}.$$ The net $(E^p_m \setminus B^p_{m, i})_{i \in
I}$ clearly tends to void set for every $m \in \mathbb{N}$ and $p
\in P$. Since $\bold m$ is a purely atomic operator-valued
measure, then according to Lemma~\ref{lemma} its semivariation is
GS-continuous, and therefore
$$ \lim_{i \in I} \hat{\bold m}^*_{p, q} (E^p_m \setminus B^p_{m, i}) =
0, \quad q \in Q, \,p \in P,\,\,\textrm{such that } q \mapsto_E
p.$$

Let $\varepsilon > 0$ be given. To every  $p \in P$ and $m \in
\mathbb{N}$ there exists an index  $j = j(m, p) \in I$, such that
for $q \mapsto _E p$, $$\hat{\bold m}_{p, q}(E^p_m \setminus
B^p_{m, i}) < \varepsilon \cdot \alpha_p \cdot \beta_m,$$ holds
for every $i \geq j(m,p)$, where $\{ \alpha_p;\, p \in P \}$ is a
summable system of positive numbers in the sense of Moore--Smith
and $\{ \beta_m; \,m \in \mathbb{N} \}$ is an absolutely
convergent series of positive numbers. Put $$F = \bigcup_{m \in
\mathbb{N}} \bigcup_{p \in P} \Bigl(E^p_m \setminus B^p_{m, j(m,
p)}\Bigr).$$ So, we have: {\setlength\arraycolsep{2pt}
\begin{eqnarray*}
\hat{\bold m}^*_{p, q}(F) & = &  \hat{\bold m}^*_{p, q}
\left(\bigcup_{m \in \mathbb{N}} \bigcup_{p \in P}\Bigl(E_m^p
\setminus B^p_{m, j(m, p)}\Bigr) \right) \\ & = & \lim_{ K = \{
p_1, \dots, p_m\}}
\hat{\bold m}_{p, q} \left( \sum_{m = 1}^{\infty} \sum_{p \in K}\Bigl( E_m^p \setminus B^p_{m, j(m, p)} \Bigr) \right) \\
& \leq & \lim_{ K = \{ p_1, \dots, p_m \} } \sum_{m = 1}^{\infty}
\sum_{p \in K} \hat{\bold m}_{p, q} \Bigl(E_m^p \setminus B^p_{m,
j(m, p)}\Bigr)
\\ & < & \varepsilon.
\end{eqnarray*}}

Let us show that the convergence of net of functions $(\bold
f_i)_{i \in I}$ is uniform on $E \setminus F$. Note that
$\bigcup_{m = 1}^{\infty} E_m^p = E$ for every $p \in P$. For a
given $\eta > 0$ choose an $m_0 \in \mathbb{N}$, such that
$\frac{1}{m} < \eta$. Then {\setlength\arraycolsep{2pt}
\begin{eqnarray*}
E \setminus F & = & E \setminus \bigcup_{m \in
\mathbb{N}} \bigcup_{p \in P}\Bigl(E_m^p \setminus B^p_{m, j(m, p)}\Bigr) \\
& = & \bigcap_{m \in \mathbb{N}} \bigcap_{p \in P} B^p_{m, j(m,
p)} \subset  B^p_{m_0, j(m_0, p)}
\end{eqnarray*}}for every $p \in P$. By definition of the set $B^p_{m_0, j(m_0, p)}$
we have that if $t \in B^p_{m_0, i}$, then
$$p(\bold f_i(t) - \bold f(t)) < \eta$$ for every
$i \geq j(m_0, p)$. So,~(\ref{net1}) implies that for every $\eta
> 0$ and $p \in P$ there exists an index $j = j(\eta, p)$, such
that for every $i \geq j(\eta, p)$, $i \in I$, there is
$$p(\bold f_i(t) - \bold f(t)) < \eta, \,\,\,t \in E \setminus B^p_{m_0, i} \supset E \setminus F,$$
i.e., the net $(\bold f_i)_{i \in I}$ converges uniformly on $E
\setminus F$. The proof is complete. \qed


\vspace{5mm}

\noindent \small{J\'an Halu\v ska, Mathematical Institute of
Slovak Academy of Science, {\it Current address:} Gre\v s\'akova
6, 040 01 Ko\v sice, Slovakia
\newline {\it E-mail address:} jhaluska@saske.sk}

\vspace{5mm}

\noindent \small{Ondrej Hutn\'ik, Institute of Mathematics,
Faculty of Science, Pavol Jozef \v Saf\'arik University in Ko\v
sice, {\it Current address:} Jesenn\'a 5, 040 01 Ko\v sice,
Slovakia,
\newline {\it E-mail address:} ondrej.hutnik@upjs.sk}

\end{document}